\documentclass{amsart}

\usepackage{amsfonts}
\usepackage{amsmath}
\usepackage{amsthm}
\usepackage{amssymb}
\usepackage[english]{babel}
\usepackage{graphics}
\usepackage{latexsym}
\usepackage{longtable} \usepackage{mathrsfs}
\usepackage{graphicx}
\usepackage{wrapfig} 
\usepackage[pdftex,colorlinks=true]{hyperref}

\newtheorem{theorem}{Theorem}
\newtheorem{corollary}[theorem]{Corollary}
\newtheorem{lemma}[theorem]{Lemma}

\newtheorem{claim}[theorem]{Claim}
\newtheorem{example}[theorem]{Example}

\theoremstyle{definition}
\newtheorem{definition}[theorem]{Definition}

\newcommand{\mL}{\mathcal{L}}

\newcommand{\mK}{\mathcal{K}}

\newcommand{\mM}{\mathcal{M}}

\newcommand{\R}{\mathbb{R}}

\newcommand{\noi}{\noindent}
\newcommand{\ms}{\medskip}

\newcommand{\al}{\alpha}
\newcommand{\be}{\beta}
\newcommand{\ga}{\gamma}
\newcommand{\de}{\delta}
\newcommand{\De}{\Delta}

\newcommand{\Om}{\Omega}

\newcommand{\larrow}{\longrightarrow}
\newcommand{\ot}{\otimes}

\newcommand{\ri}{\rightarrow}

\newcommand{\p}{\partial}
\newcommand{\sub}{\subseteq}

\newcommand{\by}{\times}
\newcommand{\rk}{\textrm{rk}}

\newcommand{\tr}{\textrm{tr}}

\newcommand{\Div}{\textrm{Div}}

\newcommand{\co}{\overline{\textrm{co}}}

\newcommand{\bt}{\begin{theorem}}\newcommand{\et}{\end{theorem}}
\newcommand{\bd}{\begin{definition}}\newcommand{\ed}{\end{definition}}
\newcommand{\bl}{\begin{lemma}}\newcommand{\el}{\end{lemma}}
\newcommand{\beq}{\begin{equation}}\newcommand{\eeq}{\end{equation}}
\newcommand{\bc}{\begin{claim}}\newcommand{\ec}{\end{claim}}
\newcommand{\bex}{\begin{example}}\newcommand{\eex}{\end{example}}
\newcommand{\bcor}{\begin{corollary}}\newcommand{\ecor}{\end{corollary}}
\newcommand{\bp}{\begin{proof}}\newcommand{\ep}{\end{proof}}

\newcommand{\BPL}{\medskip \noindent \textbf{Proof of Lemma} }
\newcommand{\BPC}{\medskip \noindent \textbf{Proof of Claim} }

\numberwithin{equation}{section}

\begin{document}

\title[Nonuniqueness in Vector-Valued Calculus of Variations in $L^\infty$...]{Nonuniqueness in Vector-Valued Calculus of Variations in $L^\infty$\\ and some Linear Elliptic Systems}

\author{\textsl{Nikos Katzourakis}}
\address{Department of Mathematics and Statistics, University of Reading, Whiteknights, PO Box 220, Reading RG6 6AX, Berkshire, UK}
\email{n.katzourakis@reading.ac.uk}


\date{}



\begin{abstract} For a Hamiltonian $H \in C^2(\R^{N \by n})$ and a map  $u:\Om \sub \R^n \!\larrow \R^N$, we consider the supremal functional
\[  \label{1} \tag{1}
E_\infty (u,\Om) \ :=\ \big\|H(Du)\big\|_{L^\infty(\Om)}  .
\]
The ``Euler-Lagrange" PDE associated to \eqref{1} is the quasilinear system
 \[  \label{2}
A_\infty u \, :=\, \Big(H_P \ot H_P +  H[H_P]^\bot \! H_{PP}\Big)(Du):D^2 u \, =  \, 0. \tag{2}
\]
\eqref{1} and \eqref{2} are the fundamental objects of vector-valued Calculus of Variations in $L^\infty$ and first arose in recent work of the author \cite{K1}. Herein we show that the Dirichlet problem for \eqref{2} admits for all $n=N\geq 2$ infinitely-many smooth solutions on the punctured ball, in the case of $H(P)=|P|^2$ for the $\infty$-Laplacian and of $H(P)= {|P|^2}{\det(P^\top\! P)^{-1/n}}$ for optimised Quasiconformal maps. Nonuniqueness for the linear degenerate elliptic system $A(x):D^2u =0$ follows as a corollary. Hence, the celebrated $L^\infty$ scalar uniqueness theory of Jensen \cite{J} has no counterpart when $N\geq 2$. The key idea in the proofs is to recast \eqref{2} as a first order differential inclusion $Du(x) \in \mK \sub \R^{n\by n}$, $x\in \Om$.
\end{abstract}

\maketitle

\section{Introduction} \label{section1}

Let $H\in C^2(\R^{N\by n})$ be  a smooth function which we call Hamiltonian, defined on the matrix space wherein the gradient $Du(x)=(D_iu_\al(x))^{\al =1,...,N}_{i =1,...,n}$ of mappings $u : \Om \sub \R^n \larrow \R^N$ is valued, $n,N\geq 2$. In this paper we consider the question of uniqueness for the Dirichlet problem of appropriately defined minimisers of vector-valued $L^\infty$ variational problems for the supremal functional
\beq \label{1.1}
E_\infty(u,\Om)\, :=\, \big\| H(Du)\big\|_{L^\infty (\Om)}.
\eeq
The 2nd order system which plays the role of ``Euler-Lagrange PDE" for \eqref{1.1} is 
\beq    \label{1.2}
A_\infty u \, :=\, \Big(H_P \ot H_P +H[H_P]^\bot H_{PP}\Big)(Du):D^2u\, =\, 0.
\eeq
System \eqref{1.2} in index form reads
\beq   \label{1.3}
\Big( H_{P_{\al i}}(Du) H_{P_{\be j}}(Du)  +  H(Du)[H_P(Du)]_{\al \ga}^\bot H_{P_{\ga i}P_{\be j}}(Du) \Big)\, D^2_{ij} u_\be\, = \, 0
\eeq
with triple summation in $i,j \in \{1,...,n\}$ and $\be \in \{1,...,N\}$ being implied. In \eqref{1.2} the notation $H_P$, $H_{PP}$ denotes derivatives of the Hamiltonian with respect to the matrix argument and $[H_P(Du(x))]^\bot$ denotes the orthogonal projection in $\R^N$ on the nullspace of the transpose of the linear map $H_P(Du(x)) : \R^n \larrow \R^N$, for $x\in \Om$. The Hessian tensor of $u$ is viewed as a map $D^2u : \Om \larrow \R^{N\by n^2}$.  The functional \eqref{1.1} and the system \eqref{1.2} are the fundamental objects of vector-valued Calculus of Variations in the $L^\infty$ norm and have first been systematically considered in the literature by the author in the recent papers \cite{K1}-\cite{K6}. The scalar case of $N=1$ has a long history, starting with Aronsson in the 60's, who first considered such $L^\infty$ variational problems in \cite{A1}-\cite{A5}. In this case of $N=1$, \eqref{1.3} simplifies to the single PDE
\beq   \label{1.4}
H_{P_{i}}(Du) H_{P_{ j}}(Du) D^2_{ij} u\, = \, 0
\eeq
which is known as the Aronsson equation, corresponding to the functional \eqref{1.1} for $H\in C^1(\R^n)$. The normal coefficient $H[H_P]^\bot H_{PP}$ of \eqref{1.2}, which has a strong geometric interpretation, in the scalar case vanishes. The special case of $H(P)=|P|^2$ of the Euclidean norm on $\R^n$ is a particularly important model and has attracted a lot of interest in the last 50 years. It reads
\beq   \label{1.5}
\De_\infty u \,:=\, D_iu\, D_ju\, D^2_{ij} u\, = \, 0
\eeq
and can be formally derived by expansion of the $p$-Laplacian $\Div(|Du|^{p-2}Du)=0$ which is the Euler-Lagrange PDE of the $p$-Dirichlet functional $\int_\Om H(Du)$ and passage to the limit as $p\ri\infty$. Noting that the functional \eqref{1.1} is non-local in the sense that minimisers over a domain $\Om$ are not minimisers over subdomains $\Om'\Subset\Om$, Aronsson remedied nonlocality by strengthening the minimality notion to \emph{Absolute Minimality}:
\beq   \label{1.6}
\Om'\Subset\Om,\ g\in W^{1,\infty}_0(\Om')\ \  \Longrightarrow\ \ E_\infty(u,\Om')\,\leq\, E_\infty(u+g,\Om').
\eeq
In the case of $H(P)=|P|^2$, Aronsson proved the equivalence between smooth Absolute Minimisers of \eqref{1.1} and classical solutions of the degenerate elliptic equation \eqref{1.5}. He also proved, by means of a Comparison Principle, \emph{uniqueness of smooth solutions to the Dirichlet problem for the $\infty$-Laplacian} in \cite{A4}, and hence uniqueness of smooth Absolute Minimisers. 

The $L^\infty$ breakthroughs of Aronsson remained dormant for some years, because the generally nonsmooth solutions of the nondivergence form nonlinear PDE \eqref{1.5} could not be rigorously interpreted until the discovery of the theory of viscosity solutions for fully nonlinear elliptic PDE by Crandall-Ishii-Lions. In the seminal paper \cite{J}, Jensen was the first to prove uniqueness of viscosity solutions for the Dirichlet problem for \eqref{1.5} and Aronsson's equivalence in the general case. His results have subsequently been extended and sharpened in major contributions by several authors in \cite{BB, BJW, CGW, JWY, CDP1, CDP2, Y, AS, ACJS}. For more details on the long history of the scalar case and a more complete list of references we refer to \cite{C} and the more recent \cite{ACJS}.

In the vector case of $N\geq 2$, the situation is much more intricate. \eqref{1.1} and \eqref{1.2} first arose in recent work of the author in \cite{K1}, wherein the PDE was formally derived in the limit of the Euler-Lagrange system of the $p$-Dirichlet functional as $p\ri \infty$ and many other PDE and variational properties where studied. The case $H(P):=|P|^2(=P:P=P_{\al i}P_{\al i})$ of the Euclidean norm  on $\R^{N\by n}$ in \eqref{1.1} is a particularly important convex example and corresponds to the \emph{$\infty$-Laplace system}:
\beq \label{1.7}
\De_\infty u \, :=\, \Big(Du \ot Du  + |Du|^2 [Du]^\bot \! \ot I \Big) : D^2 u\, = \, 0. 
\eeq
In \eqref{1.7}, $[Du(x)]^\bot$ simplifies to the orthogonal projection on the nullspace of the matrix $Du(x)^\top \! \in \R^{n\by N}$, which is the orthogonal complement of the tangent space of the map $u$. In index form, system \eqref{1.7} reads
\beq  \label{1.8}
 D_i u_\al  D_j u_\be D_{ij}^2 u_\be  +|Du|^2 [Du]_{\al \be}^\bot D^2_{ii} u_\be\ = \ 0.
\eeq
Except for the emergence of ``singular solutions'' to \eqref{1.1}, a further difficulty not present in the scalar case is that \emph{\eqref{1.1} has discontinuous coefficients} even for $C^\infty$ solutions, because there exist smooth solutions of \eqref{1.7} whose rank of the gradient is not constant: such an example on $\R^2$ is given by $u(x,y) = e^{ix}-e^{iy}$, which is $\infty$-Harmonic near the origin and has $\rk(Du)=1$ on the diagonal, but it has $\rk(Du)=2$ otherwise (\cite{K1}). In general, $\infty$-Harmonic maps present a phase separation: on each phase the dimension of the tangent space is constant and these phases are separated by \emph{interfaces} whereon the rank of $Du$ ``jumps'' and $[Du]^\bot$ is discontinuous (see \cite{K5}).  Much more complicated examples of interfaces  of smooth solutions with singularities are very recently demostrated in \cite{K4}.

The appropriate variational notion connecting \eqref{1.1} with \eqref{1.2} has been introduced in \cite{K2}, wherein the equivalence between classical solutions of \eqref{1.2} and smooth \emph{$\infty$-Minimal maps of \eqref{1.1}} has been established in the model case of $\infty$-Laplacian and under a simplifying rank restriction. The appropriate definition of minimality in the vector case is \emph{absolute minimality with respect to essentially scalar variations}, coupled by a notion of \emph{$\infty$-minimal area of $u(\Om)\sub \R^N$}. More precisely, in \cite{K2} we defined for a smooth map $u : \Om \sub \R^n \larrow \R^N$ of full rank (that is, $u$ is either an immersion or a submersion on $\Om$) the following: $u$ is called \emph{Rank-One Absolute Minimal on $\Om$} when $u$ is a local minimizer on $\Om$ with respect to essentially scalar variations $u+g\xi$ with fixed boundary values:
\beq \label{2.25}
\left.
\begin{array}{c}
\Om' \Subset \Om, \\
g\in C^1_0(\Om'),\\
\xi \in \R^N,\ |\xi|=1
\end{array}
\right\} \ \ \Longrightarrow \ \
E_\infty(u,\Om') \ \leq \ E_\infty(u+g\xi ,\Om') .
\eeq
\[
\underset{\text{Figure 1.}}{\includegraphics[scale=0.18]{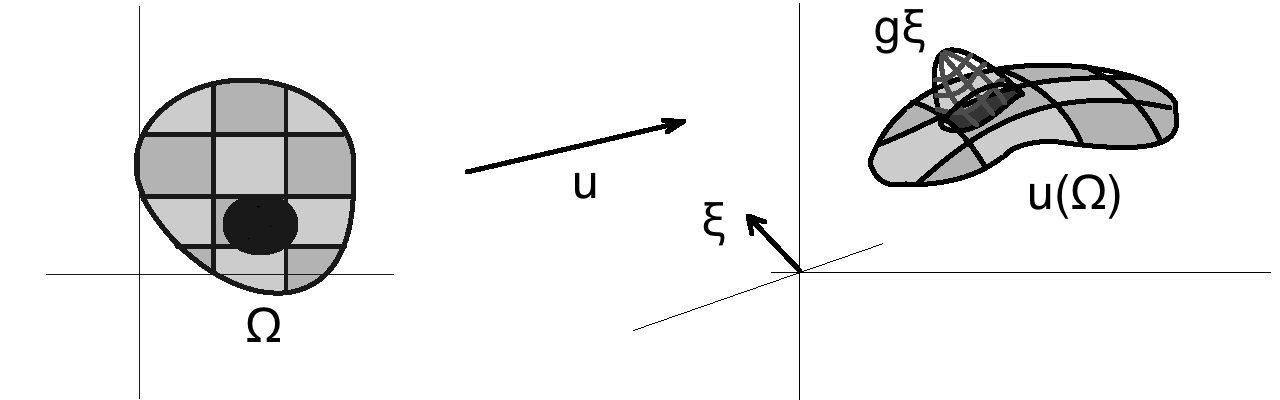}}
\]
We say that \emph{$u(\Om)$ has $\infty$-Minimal Area} when $u$ is a local minimizer on $\Om$ with respect to free variations $u+h\nu$ normal to $u(\Om)$:
\beq \label{2.26}
\left.
\begin{array}{c}
\Om' \Subset \Om, \\
 h\in C^1(\overline{\Om'}), \\
\nu \in C^1(\Om)^N,\ \nu_\al D_iu_\al=0
\end{array}
\right\} \ \ \Longrightarrow \ \
E_\infty(u,\Om')\ \leq \ E_\infty(u+h\nu,\Om').
\eeq
\[
\underset{\text{Figure 2.}}{\includegraphics[scale=0.18]{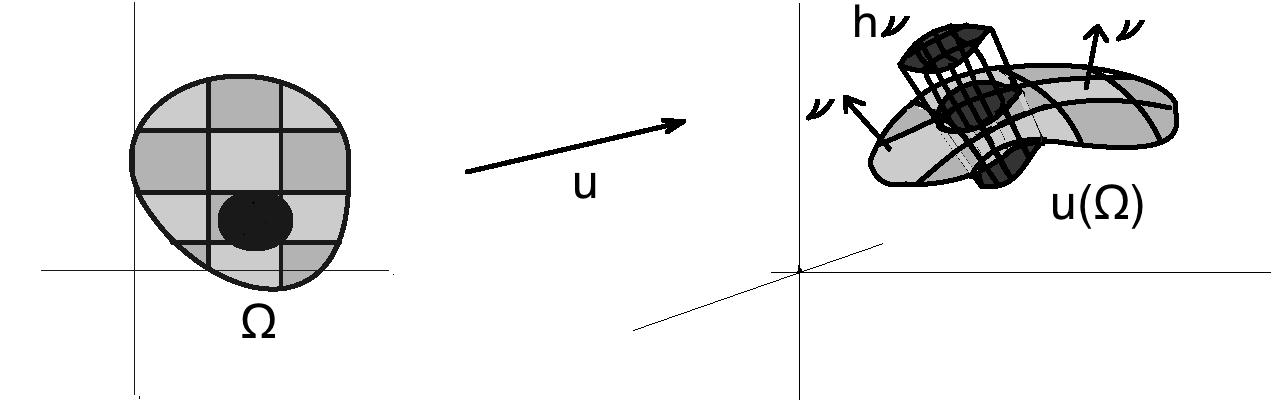}}
\]
We denote the set of $\infty$-Minimal Maps satisfying \eqref{2.25} and \eqref{2.26} by
\[ 
\mM_\infty(\Om)^N.
\]
In the scalar case of $N=1$, the set $\mM_\infty(\Om)$ of Minimal functions coincides with the set of Absolute Minimizers introduced by Aronsson. Extensions of definitions and results of \cite{K2} to a general convex Hamiltonian $H\in C^2(\R^{N\by n})$ can be established without much dificulty. The main result of \cite{K2} was 

\ms
\noi \textbf{Theorem} [$L^\infty$ variational characterisation of the $\infty$-Laplace system] \emph{Suppose $u : \Om \sub \R^n \larrow \R^N$ is a map of full rank in $C^2(\Om)^N$. Then, we have the implication}
\[
 u \in \mM_\infty(\Om) \ \ \Longrightarrow \ \ \De_\infty u=0 \text{ on }\Om ,
\]
\emph{that is, if $u$ is $\infty$-Minimal on $\Om$, then $u$ is $\infty$-Harmonic  on $\Om$. The reverse implication is true in the class of immersions.}
\ms

Extensions of the results of the papers \cite{K1,K2} to the subelliptic setting of Carnot-Carath\'eodory spaces recently appeared in \cite{K3}.

Herein we are interested in the fundamental question of uniqueness of solutions to the Dirichlet problem
for the system \eqref{1.2}:
\beq
\left\{
\begin{array}{r}
A_\infty u \, =\, 0, \ \text{ in }\Om,\ \\
u\,=\, b\ \text{ on }\p \Om,
\end{array}
\right.
\eeq
where $\Om \sub \R^n$ is a bounded connected domain and $b\in C^0(\overline{\Om})^N$. To this end, we consider two particular cases of Hamiltonians. The first one is the Euclidean norm on $\R^{N\by n}$. Surprisingly, in all dimensions $n=N\geq 2$, even smooth solutions whose restriction on the punctured ball $\Om$ equals the identity and are diffeomorphisms of $\Om$ to itself are nonunique. Accordingly, the first main result herein is

\bt[Infinitely-many solutions for $\De_\infty$] \label{th1} Let $u: \Om \sub \R^n \larrow \R^n$ be a map in $C^2(\Om)^n$, $n\geq 2$. Then, the Dirichlet problem for the $\infty$-Laplace system
\beq \label{eq2.1}
\left\{
\begin{array}{r}
\De_\infty u \ =\ 0, \ \ \text{ in } \Om,\ \, \\
u(x)\ =\ x,\  \text{ on } \p \Om,
 \end{array}
\right.
\eeq
admits infinitely-many solutions in  $C^\infty(\Om)^n$ when $\Om=\{x\in \R^n : 0<|x|<1\}$ is the punctured ball. Moreover, these solutions are diffeomorphisms from $\Om$ to itself of the form $u(x)=e^{g(|x|^2)}x$ for some $g:(0,\infty)\ri \R$.
\et

Theorem \ref{th1} is in direct contrast to the situation of the scalar case when $N=1$. Although solutions are nonunique, in \cite{K5} we showed that for $n=N=2$ they \emph{do} however satisfy a vectorial version of the Maximum Principle, well known in Calculus of Variations and Minimal surfaces  as the \emph{Convex Hull Property} (see \cite{K}):
\beq
\ \ \ \ \ u(\Om')\, \sub\, \co\big( u(\p \Om')\big), \ \ \Om' \Subset \Om.
\eeq
The failure of uniqueness for \eqref{1.7} should not be interpreted as a defect of vector-valued $L^\infty$ variational problems arising from a convex Hamiltonian, but instead as that in the absence of rigidity when $N\geq 2$ the standard scalar question may not be the appropriate one to be asked and further restrictions should be posed apart from the Dirichlet boundary condition in order to have well-posedness for the nondivergence form quasilinear system \eqref{1.7}. 

By utilising the variational characterisation of the paper \cite{K2} and Theorem \ref{th1} above, we obtain that, unlike the case of $N=1$, when $N\geq 2$ the archetypal vector-valued minimisation problem for the $L^\infty$ norm of the gradient does not have a unique solution:

\bcor[Nonuniqueness of $\infty$-Minimal Maps] \label{cor1} Let $\Om=\{x\in \R^n : 0<|x|<1\}$ be the punctured ball, $n\geq 2$. Then, the variational problem
\beq  \label{eq2.23}
\left\{ 
\begin{array}{c}
u\, \in \, \mM_\infty(\Om)^n, \ms\\
u(x)\, = \, x, \ \ x\in \p \Om,
 \end{array}
\right.
\eeq
admits infinitely-many smooth solutions $u: \Om \sub \R^n \larrow \R^n$ in $C^\infty(\Om)^n$ which are diffeomorphisms from $\Om$ to itself.
\ecor

The second Hamiltonian we consider is more complicated because it is polyconvex but not nonconvex, assuming the value $+\infty$ as well. We define the \emph{dilation function} $K : \R^{N \by n} \ri [0,+\infty]$ by
\beq \label{3.1}
K(P)\ :=\ \left\{\begin{array}{l}
           \dfrac{|P|^2}{\det\big(P^\top \! P \big)^{1/n}}\ ,  \ \ \ \  \, \rk(P)=n, \ P \in \R^{N \by n},\ms\\
          +\infty \ ,  \hspace{60pt} \rk(P)<n, \ P \in \R^{N \by n},
           \end{array}
\right.
\eeq
and we consider the respective supremal dilation functional
\beq \label{3.2}
K_\infty(u,\Om)\ :=\ \big\|K(Du)\big\|_{L^\infty(\Om)}
\eeq 
defined on (weakly differentiable) immersions $u : \Om \sub \R^n \larrow \R^N$. By utilizing the simple algebraic fact that 
\[
K(P)\geq n, \ \text{ and } \ \  K(P)=n \ \ \Longleftrightarrow  \ \ P^\top\! P = \frac{1}{n}|P|^2I
\]
it follows that $K_\infty$ measures the deviation of the immersion $u$ from being conformal on $\Om$, in which case would satisfy $Du^\top \! Du =\frac{1}{n}|Du|^2 I$. Hence, $u$ is conformal if and only if $K(Du)\equiv n$ on $\Om$. Mappings for which $K_\infty(u,\Om)<\infty$ are called \emph{Quasiconformal on $\Om$}. The case $n=N$ has a particularly long history in Geometric Analysis \cite{Ah1, Ah2, AIM, AIMO, B, G, S, T, V}. However, mostly Quasiconformal maps with distortion bounded in $L^p$ have been considered, due to the absense of a vector-valued $L^\infty$ theory. The motivation to study extremal problems for deformations of domains $\Om \ri u(\Om)$ by Quasiconformal maps is to fill the gap in between Quasiconformal maps which are too many and Conformal maps which are very few. The idea that $L^\infty$ variational mathods can be used to study ``optimised" (or ``extremal") quasiconformal maps with  \eqref{3.1} as concrete example of Hamiltonian first arose very recently in the fundamental work of Capogna and Raich \cite{CR}, almost simultaneously to the author's work in the generality of \eqref{1.1}. Actually, the results and the proofs herein are largely inspired and motivated by the ideas and the results appearing in the paper \cite{CR}. The results of \cite{CR} have been subsequently sharpened and extended by the author in \cite{K6}. The main result of \cite{K6} was the equivalence between a variational notion called \emph{``$\infty$-Minimal Quasiconformal maps"} and solutions to the system
\beq  \label{3.4}
Q_\infty u \ :=\ \Big(K_P \ot K_P +[K_P]^\bot K_{PP}\Big)(Du):D^2u\, =\, 0.
\eeq
The derivatives $K_P,K_{PP}$ of \eqref{3.1} are given by formulas \eqref{3.5} and \eqref{3.6} in Section \ref{section3}, and have been calculated in \cite{K6}. Following \cite{K6}, we call the solutions of \eqref{3.4} ``Optimal $\infty$-Quasiconformal maps". The reason that $K$ does not appear in front of $[K_P]^\bot$ as in the general case of \eqref{1.2} is that it is strictly positive and we may use a different scaling that removes it, since the two summands comprising the coefficient of $D^2u$ in \eqref{1.2} are normal to each other (see \cite{K1} and the proofs in the subsequent sections). Our second main result herein is

\bt[Infinitely-many solutions for $Q_\infty$] \label{th2} Let $u: \Om \sub \R^n \larrow \R^n$ be a map in $C^2(\Om)^n$, $n\geq 2$. Then, the Dirichlet problem for the system describing Optimal $\infty$-Quasiconformal maps
\beq \label{eq3.8}
\left\{
\begin{array}{r}
Q_\infty u \ =\ 0, \ \ \text{ in } \Om,\ \, \\
u(x)\ =\ x,\  \text{ on } \p \Om,
 \end{array}
\right.
\eeq
admits infinitely-many solutions in  $C^\infty(\Om)^n$ when $\Om=\{x\in \R^n : 0<|x|<1\}$ is the punctured ball. Moreover, these solutions are diffeomorphisms from $\Om$ to itself of the form $u(x)=|x|^\ga x$ for $\ga>-1$.
\et

In the generality of the nonconvex Hamiltonian \eqref{3.1}, uniqueness were not to be expected, at least not without extra conditions. On the other hand, Conformal maps are rigid enough so that even in the case of $n=N=2$ where the Riemann mapping theorem holds, we can not generally prescribe boundary values.

Our last auxiliary result is the proof of nonuniqueness for the Dirichlet problem for the general \emph{linear} homogeneous degenerate elliptic 2nd order system in nondivergence form:
\beq \label{4.1}
A_{\al i \be j}(x)D^2_{ij}u_\be (x) \, =\, 0 .
\eeq
Here $A_{\al i \be j}$ are the components of a 4th order tensor which is smooth, symmetric and convex as a linear map on $\R^{N \by n}$:
\beq \label{4.2}
A_{\al i \be j} \in C^{\infty}(\Om), \ \ A_{\al i \be j}=A_{\be j \al i }, 
\eeq
\beq  \label{4.3}
0\ \leq\ A_{\al i \be j}(x)P_{\al i}P_{\be j}\ \leq \ c\, |P|^2, 
\eeq
for some $c>0$, $P\in \R^{N \by n}$, $ x\in \Om$. In corollary \ref{pr1} below we will show that even when $n=N$ and $A_{\al i \be j} $ satisifes the partial ellipticity
\beq \label{4.5}
A_{\al i \be j}(x)P_{\al i}P_{\be j}\ \geq \ \Big(\left( [x]^\bot\! + (1-\mu)[x]^\top\right):P\Big)^2,
\eeq
solutions to \eqref{4.1} are nonunique. Here partial ellipticity is meant in the sense that we have convexity along certain lines of $\R^{n\by n}$ generated by the linear combination of the projections $[x]^\top$ and $[x]^\bot$, where $[x]^\top$ and $[x]^\bot$ are given by
\beq \label{eq2.14}
I\ = \ \frac{x \ot x}{|x|^2}\ +\ \left(I\, -\, \frac{x \ot x}{|x|^2}\right)\ =: \ [x]^\top +\ [x]^\bot.
\eeq
The projections refered to are orthogonal projections with respect to the Euclidean norm of $\R^{n\by n}$.

\begin{corollary}[Nonuniqueness for linear nondivergence elliptic systems] \label{pr1}
 Let $\Om=\{x\in \R^n :0<|x|<1\}$ be the punctured ball, fix $\mu >1$ and let $A_{\al i \be j}$ be the components of a 4th order tensor given by
\beq \label{4.4}
A_{\al i \be j}(x)\, :=\, \left(\de_{\al i}-\mu \frac{x_\al x_i}{|x|^2}\right)\left(\de_{\be j}- \mu\frac{x_\be x_j}{|x|^2}\right)
\eeq
which satisfies \eqref{4.2}, \eqref{4.3} and \eqref{4.5}. Then, for any $\mu>1$ and $n\geq 2$, the Dirichlet problem
\beq \label{4.6}
\left\{
\begin{array}{r}
A_{\al i \be j}(x)D^2_{ij}u_\be (x) \, =\, 0,\ \ x\in \Om, \ \\
u(x)\,=\, x,\ x\in \p \Om,
\end{array}
\right.
\eeq
admits two smooth solutions $u^\mu,u^* \in C^\infty(\Om)^n$ which are diffeomorphisms of $\Om$ to itself and are given by
\beq \label{4.7}
u^\mu(x)\, :=\, |x|^{\frac{n-\mu}{\mu -1}}x,\ \ \ u^*(x)\,=\, x.
\eeq
\end{corollary}

The above counterexample to uniqueness, has, not surprisingly, a vague resemblance to the famous counterexample to regularity of vector-valued minimisers of De Giorgi (e.g. see \cite{Gi}).

The idea of the proof in Theorems \ref{th1} and \ref{th2} is to obtain a certain class of solutions to \eqref{1.2} by constructing \emph{smooth} solutions of first order differential inclusions of the form
\beq
Du(x) \, \in \, \mK \sub \R^{N\by n},\ \ x\in \Om.
\eeq
There exist well known methods which allow to construct $W^{1,\infty}$ and $C^1$ solutions to such differential inclusions (Gromov's Convex Integration and its extension by Muller-Sverak \cite{MS}, Baire's Category method of Dacorogna-Marcellini \cite{DM}) but none of these general methods provide smooth solutions. The idea in the proof of Corollary \ref{pr1} is a simple argument of freezing the nonlinear coefficients of $D^2u$ in \eqref{1.2}.

We conclude this long introduction by noting that the only vectorial result we are aware of and strongly relates to the material herein except for the paper \cite{CR} is the paper by Sheffield and Smart \cite{SS}. Thereis they use the nonsmooth operator norm on $\R^{N \by n}$ and derive a singular variant of \eqref{1.7} which governs \emph{vector-valued optimal Lipschitz extensions}.

\section{Nonuniqueness for the $\infty$-Laplacian.} \label{section2}

This section is devoted to the proof of Theorem \ref{th1}. The key idea is that the class of smooth solutions $u:\Om \sub \R^n \larrow \R^n$ to the $\infty$-Laplace system which are local diffeomorphisms coincides with the set of smooth solutions of a \emph{first order differential inclusion}. In Lemma \ref{l2a} below we establish this equivalence and subsequently in  Lemma \ref{l2}  we construct \emph{smooth} solutions to this first order differential inclusion.

\bl[The $\infty$-Laplacian as Differential Inclusion] \label{l2a} Let  $u:\Om \sub \R^n \larrow \R^n$ be a map in $C^2(\Om)^n$, suppose $\Om$ is connected and for $a\geq0$ consider the set of matrices
\beq  \label{eq2.2}
\mL_a\, :=\, \Big\{ A \in \R^{n \by n} \, :\,  |A|=\sqrt{a},\  \det(A)\neq 0 \Big\} .
\eeq
Consider the following statements: 

\ms

\noi (a) There exists an $a\geq 0$ such that $u$ solves the differential inclusion
\beq  \label{eq2.3a}
Du(x)\, \in \,\mL_a,\ \  x\in \Om,
\eeq
classically on $\Om$.

\ms
\noi (b) $u$ solves $\De_\infty u =0$ on $\Om$.
\ms

Then, $(a)$ implies $(b)$ while the reverse implication is true when $u$ in addition is a local diffeomorphism on $\Om$.

\el

\BPL \ref{l2a}. Let $u:\Om \sub \R^n \larrow \R^n$ be in $C^2(\Om)^n$ and fix $a\geq 0$. We begin by rewriting the differential inclusion \eqref{eq2.3a} with $\mL_a$ given by \eqref{eq2.2} as a 1st order vectorial Eikonal PDE coupled by a rank constraint:
\beq \label{eq2.8a}
\left\{ 
\begin{array}{r}
\big|Du(x)\big|^2 =\, a, \ \  x\in\Om,\ \ \ \ \ \ \ms\\
\rk\big(Du(x)\big)\, =\, n\ ,  \ \  x\in\Om.\ \ \ \ \ \ \\
\end{array}
\right.
\eeq
We now rewrite the $\infty$-Laplace operator \eqref{1.7} as
\beq 
\De_\infty u\ =\ Du D\Big(\frac{1}{2}|Du|^2\Big)\ +\ |Du|^2[Du]^\bot\De u.
\eeq
Since $[Du]^\bot Du=0$, the two operators composing $\De_\infty$ are normal to each other and hence the system $\De_\infty u =0$ decouples to the pair of systems
\beq \label{eq2.5a}
\left\{ 
\begin{array}{r}
Du(x) D\Big(\dfrac{1}{2}|Du|^2\Big)(x)\, =\, 0, \ \ \  x\in\Om, \ms\\
\big|Du(x)\big|^2[Du(x)]^\bot\De u(x)\, =\, 0, \ \ \ x\in\Om. \\
\end{array}
\right.
\eeq
Suppose first $u$ solves \eqref{eq2.8a}. Then, $|Du|^2$ is constant on $\Om$ and since the rank of $Du$ equals $n$ on $\Om$,  we have that $[Du(x)]^\bot=0$ because the latter equals the orthogonal projection on the nullspace of the matrix $Du(x)^\top \! \in \R^{n \by n}$. Hence, $u$ satisifes \eqref{eq2.5a} and hence is $\infty$-Harmonic on $\Om$.

Conversely, suppose $\De_\infty u=0$ and moreover that $u$ is a local diffeomorphism on $\Om$. Then, $\rk(Du)=n$ on $\Om$, which gives that $(Du(x))^{-1}$ exists for all $x\in \Om$. Consequently, by \eqref{eq2.5a},
\beq
D\big(|Du|^2\big)\, =\, (Du)^{-1} DuD\big(|Du|^2\big)\, = \, 0.
\eeq
Since $\Om$ is connected, it follows that $|Du|$ is constant throughout $\Om$. Hence, there exists an $a\geq 0$ such that $|Du|^2=a$. Hence, $u$ satisfies \eqref{eq2.8a} and as a result solves the differential inclusion \eqref{eq2.3a}.                                 \qed

\ms

\bl[Existence for the Differential Inclusion] \label{l2} Let $\Om=\{x\in \R^n : 0<|x|<1\}$ be the punctured ball and for $a\geq 0$ consider again the set of matrices \eqref{eq2.2} and the differential inclusion \eqref{eq2.3a}. 

Then, for any $a>n$, the Dirichlet problem
\beq  \label{eq2.3}
\left\{ 
\begin{array}{r}
Dv(x)\, \in \,\mathcal{L}_a,\ \  x\in \Om,\ \, \\
v(x)\, =\, x,\ \ \, x\in \p \Om,\
 \end{array}
\right.
\eeq
admits a unique solution $v_a : \overline{\Om} \larrow \overline{\Om}$ in $C^\infty(\Om)^n$ of the form $v_a(x)=e^{g_a(|x|^2)}x$ for some increasing and concave $g_a \in C^\infty(0,\infty)$ with $g_a(0^+)=-\infty$ and $g_a(1)=0$. In addition, $v_a$ is a diffeomorphism on $\Om$.
\el

\BPL \ref{l2}. We begin by rewriting the Dirichlet problem \eqref{eq2.3} as the Dirichlet problem for the vectorial Eikonal PDE with rank constraint:
\beq  \label{eq2.4}
\left\{ 
\begin{array}{r}
\big|Dv\big|^2 =\, a, \ \, \ \ \ \ \text{ in }\big\{ x\in \R^n : 0<|x|<1 \big\},\ \ \ \ \ \ \ms\\
\rk(Dv)\, =\, n ,\  \ \ \ \  \text{ in }\big\{ x\in \R^n : 0<|x|<1 \big\},\ \ \ \ \ \ \ms\\
\ \ \ v(x)\, =\, x,  \ \ \ \  \text{ on }\big\{ x\in \R^n : |x|=1 \text{ or } x=0 \big\}.
 \end{array}
\right.
\eeq
We look for solutions $v : \Om=\{0<|x|<1\} \sub \R^n \larrow \R^n$ to \eqref{eq2.4} of the special form 
\beq \label{eq2.5}
v(x)\ =\ e^{g(|x|^2)}x 
\eeq
for some $g : (0,\infty)\sub \R \larrow \R$. By differentiating \eqref{eq2.5} we have $D_iu_j(x) = e^{g(|x|^2)} \big(\de_{ij} + 2 g'(|x|^2)x_ix_j\big)$, and hence
\beq \label{eq2.6}
Dv(x)\ =\ e^{g(|x|^2)}\Big(I\, +\, 2g'(|x|^2) x\ot x\Big).
\eeq
Thus, 
\beq \label{eq2.7}
Dv(x)^\top \! Dv(x)\ =\ e^{2g(|x|^2)}\Big[I\, +\, 4\Big(g'(|x|^2) \, +\, \big(g'(|x|^2)\big)^2 |x|^2\Big) x\ot x \Big]
\eeq
and by utilizing that $\tr(I)=n$, \eqref{eq2.7} gives
\beq \label{eq2.8}
\big|Dv(x)\big|^2\ =\ e^{2g(|x|^2)}\Big[n\, +\, 4g'(|x|^2)|x|^2 \,+\, 4\big(g'(|x|^2)\big)^2 |x|^4 \Big].
\eeq 
Moreover, by the identity \eqref{eq2.14}, for $x\neq 0$, we rewrite \eqref{eq2.6} as the following linear combination of complementary projections:
\beq \label{eq2.10}
Dv(x)\ =\ e^{g(|x|^2)}\Big[\Big(1\, +\, 2g'(|x|^2)|x|^2\Big) [x]^\top +\ [x]^\bot\Big],
\eeq
for $x\in \Om$. In view of \eqref{eq2.10} and \eqref{eq2.4}, we have that
\beq \label{eq2.11}
\rk(Dv)\, =\, n \text{ in }\Om \ \ \Leftrightarrow \ \  g'(t)\, \neq \, -\frac{1}{2t},\ 0<t<1.
\eeq
On the other hand, in view of \eqref{eq2.8} and \eqref{eq2.4}, we have that $|Dv|^2=a$ on $\Om$ if and only if
\beq
n\, +\, 4g'(|x|^2)|x|^2 \,+\, 4\big(g'(|x|^2)\big)^2 |x|^4 \ = \ ae^{-2g(|x|^2)},
\eeq
for $x\in \Om$, which means that
\beq  \label{eq2.13}
\text{ $|Dv|^2=a$ in $\Om$ } \ \Leftrightarrow \ \ \big(g'(t)t\big)^2 \, +\, g'(t)t \, +\, \frac{1}{4}\Big(n-ae^{-2g(t)}\Big) = 0, \ 0<t<1.
\eeq
The next claim completes the proof of the lemma.

\bc \label{cl1} The constrained first order boundary-value problem 
\beq  \label{eq2.14a}
\left\{ 
\begin{array}{c}
g'(t) \, =\, \dfrac{1}{2t}\Big(\sqrt{ae^{-2g(t)}-(n-1)}\, - \, 1 \Big), \ \ 0<t<1,\ms\\
g'(t)\, >\, 0,  \hspace{140pt} 0<t<1,\ms\\
g(0^+)\, =-\infty,\ \ \ \ g(1)\, =\, 0,\ \ \ \ 
 \end{array}
\right.
\eeq
has for every $a>n$ a unique smooth concave solution $g_a \in C^\infty(0,1)$.
\ec

\BPC \ref{cl1}. Fix $a>n$ and consider the initial value problem
\beq  \label{eq2.15}
\left\{ 
\begin{array}{l}
g'(t)\, =\, F\big(t,g(t)\big), \ms\\
\, g(1)\, =\, 0,
 \end{array}
\right.
\eeq
for
\beq
F(t,y)\, :=\, \dfrac{1}{2t}\left(\sqrt{ae^{-2y}-(n-1)}\, - \, 1 \right).
\eeq
It is easy to verify that the function $F$ is well-defined, strictly positive and smooth on the domain
\beq
B\, :=\, (0,+\infty) \by \left(\! -\infty,\ln\sqrt{\frac{a}{n}}\right) \, \sub \, \R^2.
\eeq
Since $F\in C^\infty(B)$ and $F>0$, there exists a unique strictly increasing  local solution $g_a \in C^\infty(t^-,t^+)$ to the ODE problem \eqref{eq2.15} for some $t^-<1<t^+$. By maximality, $g_a$ can be extended to a solution $g_a \in C^\infty(0,t^+)$ of  \eqref{eq2.15} which remains strictly increasing on its domain of definition and satisfies $g_a(0^+)=-\infty$. Moreover, by differentiating \eqref{eq2.14a} we have
\begin{align}
g''(t)\,&=\, -\dfrac{1}{2t^2}\left(\sqrt{ae^{-2g(t)}-(n-1)}\,- \, 1 \right)\, \nonumber\\
&\ \ \ \ \ \ -\,\frac{1}{4t}\frac{ag'(t)e^{-2g(t)}}{2t \sqrt{ae^{-2g(t)}-(n-1)}} \\
&=\, -\frac{g'(t)}{t}\left( 1 +\frac{ae^{-2g(t)}}{2\sqrt{ae^{-2g(t)}-(n-1)}}\right),\nonumber
\end{align}
for $0<t<1$, and as a result $g_a$ is concave since $g_a'>0$.
\[
\underset{\text{Figure 3.}}{\includegraphics[scale=0.24]{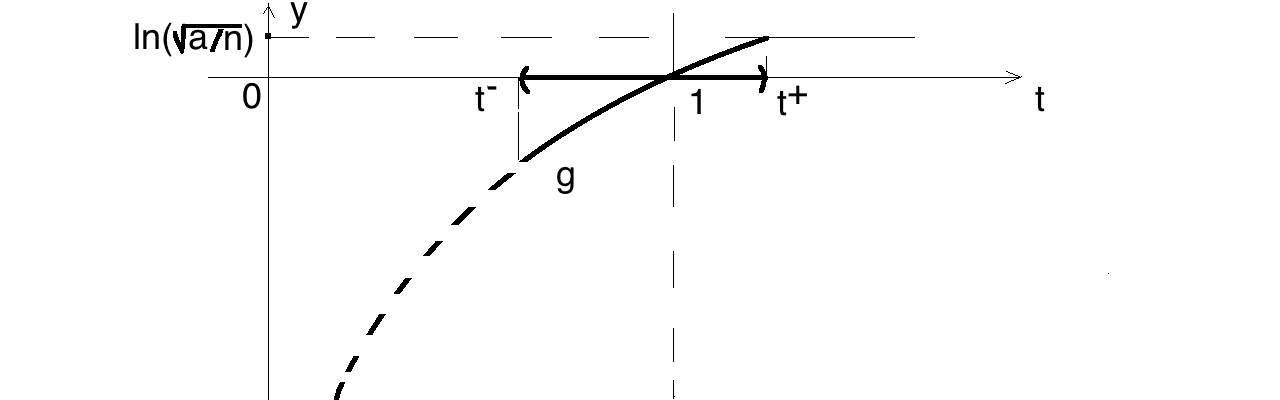}}
\]
Hence, $g_a$ solves the constrained boundary-value problem \eqref{eq2.14a} and as such the claim follows. \qed  

\ms

We may now complete the proof of Lemma \ref{l2} by observing that the mapping 
\beq \label{eq2.18}
v_a(x)\ :=\ e^{g_a(|x|^2)}x
\eeq
where $g_a$ is the solution constructed in Claim \ref{cl1}, solves the differential inclusion \eqref{eq2.3}. Indeed, we have that $v_a \in C^\infty(\Om)$, and since $g_a(0^+)=-\infty$ and $g_\al(1)=0$, we obtain that $v_a$ of \eqref{eq2.18} satisfies the desired boundary conditions because $v_a(0)=e^{-\infty}0=0$ and for $|x|=1$ we have $v_a(x)=e^{g(1)}x=x$. Moreover, in view of the equivalences \eqref{eq2.11} and \eqref{eq2.13}and the observation that the solution $g_a$ of problem \eqref{eq2.14a} satisfies the right hand side of \eqref{eq2.13}, we deduce that the map $v_a$ solves \eqref{eq2.4} and hence \eqref{eq2.3}. Finally, by utilizing that $v_a$ is bijective when restricted on radii of the unit ball, it easily follows that the local diffeomorphism $v_a :\Om \sub \R^n \larrow \R^n$ actually is a diffeomorphism of $\Om$ to itself. The lemma follows.
\qed

\ms

The proof of Theorem \ref{th1} is now complete.     
\ms

\section{Nonuniqueness of Optimal $\infty$-Quasiconformal Maps.} \label{section3}

This section is devoted to the proof of Theorem \ref{th2}. We begin by recalling from \cite{K6} that the derivatives of $K$ at the neighbourhood of a point $P \in \R^{N\by n}$ with $\rk(P)=n$ are given by
\beq \label{3.5}
K_{P_{\al i}} (P)\ =\ 2P_{\al m} \Bigg( \frac{ \de_{mi}-\frac{1}{n}|P|^2\big(P^\top \! P \big)_{mi}^{-1} }{\det\big(P^\top \! P \big)^{1/n}} \Bigg),
\eeq
and
\begin{align} \label{3.6}
K_{P_{\al i}P_{\be j}} (P)\ =\ & 2\de_{\al \be}  \Bigg( \frac{\big(P^\top \! P \big)_{ik}^{-1} \big( S(P^\top\! P)\big)_{kj}  }{\det\big(P^\top \! P \big)^{1/n}} \Bigg) \nonumber\\
& +\ 2P_{\al m}P_{\be l }   \Bigg( \frac{ \big(P^\top \! P \big)_{ik}^{-1}  E_{kjlm}  }{\det\big(P^\top \! P \big)^{1/n}}\Bigg) \ + \ O_{\al i \be j}(P).
\end{align}
Here $O_{\al i \be j}(P)$ is a tensor of the form $K_{P_{\al m}}(P)A_{m\be i j}(P)$ which is annihilated by $[K_P]_{\ga \al}^\bot$ (that is $[K_P]^\bot O=0$) and does not appear in the system \eqref{3.4}. $E$ denotes a constant 4-tensor $E_{kjlm}$ given by
\beq  \label{3.7}
E_{kjlm}\ := \ \de_{ml}\de_{jk} \, +\, \de_{mj}\de_{kl} \, -\, \frac{2}{n}\de_{mk}\de_{jl} ,
\eeq
and for $A \in \R^{n \by n}$, $S(A)$ is the \emph{Ahlfors operator} of $A$, given by
\beq \label{3.3}
S(A)\ :=\ \frac{A+A^\top}{2}\ -\ \frac{\tr(A)}{n}I
\eeq
and is a symmetric traceless matrix. The method followed to prove Theorem \ref{th2} is analogous to that followed in Theorem \ref{th1}. We begin by showing that a special class of smooth solutions to system \eqref{3.4} can be obtained by means of a differential inclusion.

\bl[Optimal $\infty$-Quasiconformal maps via a Differential Inclusion] \label{l3a} Let  $u:\Om \sub \R^n \larrow \R^n$ be a local diffeomorphism in $C^2(\Om)^n$, suppose $\Om$ is connected and for $a\geq n$ consider the set of matrices
\beq  \label{3.9}
\mathcal{K}_a\, :=\, \left\{ A \in \R^{n \by n}: \frac{\tr(A^\top \! A )\ \ }{\det(A^\top \! A )^{1/n}}=a,\ \det\big(S(A^\top\! A)\big)>0\right\} 
\eeq
where $S$ the Ahlfors operator. Consider the following statements: 

\ms

\noi (a) There exists an $a\geq n$ such that $u$ solves the differential inclusion
\beq  \label{3.10}
Du(x)\, \in \,\mathcal{K}_a,\ \  x\in \Om,
\eeq
classically on $\Om$.

\ms
\noi (b) $u$ solves $Q_\infty u =0$ on $\Om$.
\ms

Then, $(a)$ implies $(b)$ while the converse is true when in addition the Ahlfors operator $S(Du^\top\! Du)$ of the induced Riemannian metric $Du^\top\! Du$ has rank equal to $n$, everywhere on $\Om$.

\el

\BPL \ref{l3a}. The proof follows similar lines to that of Lemma \ref{l2a} and hence we merely schetch it briefly. Let $u:\Om \sub \R^n \larrow \R^n$ be in $C^2(\Om)^n$ and fix $a\geq n$. In view of \eqref{3.1}, we may rewrite the differential inclusion \eqref{3.10} with $\mK_a$ given by \eqref{3.9} as a vectorial Hamilton-Jacobi PDE with coefficients given by \eqref{3.1}, coupled by two rank constraints, on  $Du$ and on the Ahlfors operator of $Du^\top\! Du$:
\beq \label{3.11}
\left\{ 
\begin{array}{r}
K\big(Du(x)\big)\, =\, a, \ \  \ x\in\Om, \ \ \ \ \ \ms\\
\rk \big(Du(x)\big)\, =\, n\ ,  \ \  x\in\Om,\ \ \ \ \ \ms\\
\rk \big(S(Du^\top\! Du)(x)\big)\, =\, n\ ,  \ \  x\in\Om.\ \ \ \ \ \\
\end{array}
\right.
\eeq
We also rewrite the system \eqref{3.4} as
\beq  \label{3.12}
Q_\infty u\ =\ K_P(Du) D\big(K(Du)\big)\ +\ \big([K_P]^\bot K_{PP}\big)(Du):D^2 u.
\eeq
Since $[K_P]^\bot K_P=0$, the two operators composing $Q_\infty$ are normal to each other and hence \eqref{3.4} decouples to the pair of systems
\beq \label{3.13}
\left\{ 
\begin{array}{r}
K_P(Du) D\big(K(Du)\big) \, =\, 0, \ \ \  x\in\Om, \ms\\
\big([K_P]^\bot K_{PP}\big)(Du):D^2 u \, =\, 0, \ \ \ x\in\Om. \\
\end{array}
\right.
\eeq
By \eqref{3.5} and \eqref{3.3} we have
\beq
K_P(Du)\, =\, 2Du \big(Du^\top\! Du\big)^{-1}\frac{S(Du^\top\! Du)}{\det (Du^\top\! Du )^{1/n}}
\eeq
and hence the 2nd and 3rd equations of \eqref{3.11} are equivalent to 
\beq
\rk\big(K_P(Du)\big)\, =\, n,\ \ \text{ on }\Om. 
\eeq
Hence, $[K_P(Du(x))]^\bot=0$ because it equals the orthogonal projection on the nullspace of $K_P(Du(x))^\top \!\in \R^{n \by n}$. Moreover, the latter matrix is invertible, which by connectivity of $\Om$ implies $K_P(Du) D\big(K(Du)\big)=0$ if and only $K(Du)\equiv a$ for some $a\geq n$.                            \qed

With the next result we complete the proof of Theorem \ref{th2} by showing existence of infinitely-many solutions to the Dirichlet problem for the differential inclusion \eqref{3.10}, \eqref{3.9}. Material related to Lemma \ref{l3} below appears in the paper \cite{CR} as well.

\bl[Existence for the Differential Inclusion] \label{l3} Let $\Om=\{x\in \R^n : 0<|x|<1\}$ be the punctured ball and for $a> n\geq 2$ consider the set of matrices \eqref{3.9} and the differential inclusion \eqref{3.10}. 

Then, for any $\ga>-1$ and for
\beq \label{3.15}
a(\ga)\ :=\ \frac{n\, +\, \ga^2\, +\, 2\ga}{(1\, +\, \ga)^{n/2}},
\eeq
the Dirichlet problem
\beq  \label{3.16}
\left\{ 
\begin{array}{r}
Dv^\ga(x)\, \in \,\mathcal{K}_{a(\ga)},\ \  x\in \Om,\ \, \\
v^\ga(x)\, =\, x,\ \ \ \ \ \, x\in \p \Om, \
 \end{array}
\right.
\eeq
admits a solution $v^\ga : \overline{\Om} \larrow \overline{\Om}$ in $C^\infty(\Om)^n$ of the form $v^\ga (x)=|x|^\ga x$, which is a diffeomorphism on $\Om$.
\el

\BPL \ref{l3}.  Define $v^\ga (x):=|x|^\ga x$. For $x \in \Om$, we readily have
\begin{align}
Dv^\ga(x)\ &= \ |x|^\ga\left(I\, +\, \ga \dfrac{x \ot x}{|x|^2}\right), \label{3.17}\\
\big(D{v^\ga}^\top \! Dv^\ga\big)(x) \ &= \ |x|^{2\ga}\left(I\, +\, (\ga^2+2\ga) \dfrac{x \ot x}{|x|^2}\right). \label{3.18}
\end{align}
Hence, \eqref{3.18} gives
\begin{align} \label{3.19}
\big|Dv^\ga(x)\big|^2\ = \ |x|^{2\ga} (n\, +\, \ga^2+2\ga )
\end{align}
and by utilizing Sylvester's determinant theorem, we also have
\begin{align} \label{3.20}
\det\big(D{v^\ga}^\top \! Dv^\ga\big)(x) \ &= \ |x|^{2n\ga} \left[1+ (\ga^2+2\ga) \dfrac{x}{|x|}^\top\! \dfrac{x}{|x|}\right] \nonumber\\
& = \ |x|^{2n\ga} (1 +\ga)^2 .
\end{align}
By \eqref{3.19}, \eqref{3.20} and \eqref{3.1}, we obtain that $v^\ga$ has constant dilation $K(Dv^\ga)$ throughout $\Om$ and equal to $a(\ga)$, the constant given by \eqref{3.15}. By recalling the split of the identity of $\R^n$ given by \eqref{eq2.14}, \eqref{3.17} implies
\beq
Dv^\ga(x)\ = \ |x|^\ga\Big((1 +\ga) [x]^\top + \, [x]^\bot\Big)
\eeq
Since $\ga>-1$, it follows that $\rk(Dv^\ga) = n$ in $\Om$. By \eqref{3.17}, \eqref{3.3} and \eqref{eq2.14}, we calculate
\begin{align}
S\big(D{v^\ga}^\top \! Dv^\ga\big)(x) \ &= \ |x|^{2\ga}  \left(I\, +\, (\ga^2+2\ga) \dfrac{x \ot x}{|x|^2} \right)\, - \, \frac{1}{n}|x|^{2\ga} \big(n\, +\, \ga^2+2\ga \big) I \nonumber\\
&= \ (\ga^2+2\ga) |x|^{2\ga} \left(  \dfrac{x \ot x}{|x|^2}\ - \ \frac{1}{n} I \right)  \label{3.22}\\
&= \ (\ga^2+2\ga) |x|^{2\ga} \left[ \Big(1  - \frac{1}{n}\Big) [x]^\top - \  \frac{1}{n} [x]^\bot\right]. \nonumber
\end{align}
By \eqref{3.22}, we obtain that $\rk\big(D{v^\ga}^\top \! Dv^\ga\big)(x) = n$ in $\Om$. Hence, it follows that the map $v^\ga$ solves the 1st order PDE system \eqref{3.12}, which is equivalent to the differential inclusion \eqref{3.10} with $\mK_{a}$ given by \eqref{3.9} with $a=a(\ga)$ given by \eqref{3.15}. The lemma follows.                   \qed
 \ms

The proof of Theorem \ref{th2} is complete.

\section{Nonuniqueness for Linear Nondivergence Degenerate\\ Elliptic PDE Systems.} \label{section4}

In this brief section we prove Corollary \ref{pr1}.  The idea is to write the linear system \eqref{4.1} in the form or the quasilinear system \eqref{3.4} by freezing the coefficient $K_P(Du)\ot K_P(Du)$ at the map $u^\mu$ of \eqref{4.7}. By \eqref{3.5}, (since $n=N$) we first calculate
 \begin{align} \label{4.8}
K_P(P)\ &=\ 2P \big(P^\top \! P \big)^{-1} \Bigg( \frac{ P^\top \! P -\frac{1}{n}|P|^2I }{\det (P^\top \! P )^{1/n}} \Bigg)  \nonumber \\
& = \ 2P^{\top,-1}\Bigg( \frac{ S(P^\top \! P )}{\det(P^\top \! P)^{1/n}} \Bigg).
\end{align}
We now set
\beq \label{4.9}
\ga(\mu)\, :=\, \frac{n-\mu}{\mu-1}
\eeq
and then by \eqref{4.7}, \eqref{4.8} and \eqref{3.17} we have 
\beq \label{4.9a}
 u^\mu \, =\, v^{\ga(\mu)}
\eeq
(where $v^\ga(x)=|x|^\ga x$) and hence
\beq  \label{4.10}
Du^\mu(x)\ = \ |x|^{\ga(\mu)}\left(I\, +\, \ga(\mu) \dfrac{x \ot x}{|x|^2}\right).
\eeq
Hence, by \eqref{4.8}, \eqref{4.10} and \eqref{3.20}, \eqref{3.22}, we have
 \begin{align} \label{4.11}
K_P(Du^\mu)\  &= \ \frac{2 \ga(\mu) (2+ \ga(\mu) )|x|^{2 \ga(\mu) }}{ |x|^{2 \ga(\mu)} (1+ \ga(\mu) )^{2/n} }\, \cdot  \nonumber \\
&\ \ \ \ \ \ \
\cdot \left\{|x|^{\ga(\mu)}\left(I\, +\, \ga(\mu) \dfrac{x \ot x}{|x|^2} \right)\right\}^{-1}\left( \frac{x\ot x}{|x|^2} -\frac{1}{n}I \right)    \\
 &= \ -\frac{2 \ga(\mu) ( \ga(\mu) +2)}{n(1+ \ga(\mu) )^{2/n} } |x|^{-\ga(\mu)}  \left(I\, +\, \ga(\mu) \dfrac{x \ot x}{|x|^2} \right)^{-1}\left( I-n\frac{x\ot x}{|x|^2}\right). \nonumber
\end{align}
By utilizing the identity 
\beq \label{4.12}
 \left(I\, +\, \ga(\mu) \dfrac{x \ot x}{|x|^2} \right)^{-1}=\  I\, -\, \left(\frac{\ga(\mu)}{\ga(\mu)+1} \right)\dfrac{x \ot x}{|x|^2} 
\eeq
\eqref{4.11} together with \eqref{4.9} give
 \begin{align} \label{4.13}
K_P(Du^\mu)\ &= \ -\frac{2\ga(\mu)(2+\ga(\mu))}{n (1+\ga(\mu))^{2/n} } |x|^{-\ga(\mu)}  \left[ I- \left(\frac{\ga(\mu)+n}{\ga(\mu)+1} \right) \frac{x\ot x}{|x|^2}\right] \nonumber\\
&= \ -\frac{2\ga(\mu)(2+\ga(\mu))}{ n(1+\ga(\mu))^{2/n} } |x|^{-\ga(\mu)}  \left( I- \mu \frac{x\ot x}{|x|^2}\right).
\end{align}
By \eqref{4.13} and \eqref{4.4} we have
 \begin{align} \label{4.14}
 K_{P_{\al i}}(Du^\mu) K_{P_{\be j}}(Du^\mu)\ &= \ \frac{4\ga(\mu)^2(2+\ga(\mu))^2}{ n^2(1+\ga(\mu))^{4/n} } |x|^{-2\ga(\mu)} \cdot \nonumber \\
&\ \ \ \ \ \cdot \left( \de_{\al i}- \mu \frac{x_\al x_i}{|x|^2}\right)  \left( \de_{\be j}- \mu \frac{x_\be x_j}{|x|^2}\right) \nonumber \\
&=\  \frac{4\ga(\mu)^2(2+\ga(\mu))^2}{n^2 (1+\ga(\mu))^{4/n} } |x|^{-2\ga(\mu)} A_{\al i \be j}(x).
\end{align}
Hence, by Lemma \ref{l3} and \eqref{4.7}, \eqref{4.9a}, we have that $u^\mu$ is a solution to \eqref{4.1} on $\Om$, since $\ga(\mu)>-1$ and
 \begin{align} \label{4.14}
A_{\al i \be j}(x)D_{ij}^2u^\mu_{\be}(x)\ &=\ \frac{ n^2(1+\ga(\mu))^{4/n}|x|^{-2\ga(\mu)} } {4\ga(\mu)^2(2+\ga(\mu))^2} K_{P_{\al i}}(Du^\mu) K_{P_{\be j}}(Du^\mu) D_{ij}^2u^\mu_{\be} \nonumber\\
&=\ 0.
\end{align}
Moreover, $u^*(x)=x$ is also a solution to \eqref{4.1}, since $D^2u^* \equiv 0$. The rest of the claims follow easily by arguing as in Theorem \ref{th2}. The corollary follows.             \qed

\ms

\noi \textbf{Acknowledgement.} The author wishes to thank Y. Yu and J. Manfredi for their suggestions and comments as well as Ch. Wang for his selfless share of expertise of the subject and the scientific discussions during the author's visit to the Department of Mathematics at Kentucky. He is also indebted to L. Capogna for the scientific discussions about extremal problems in Geometric Analysis.

\bibliographystyle{amsplain}

\end{document}